\documentclass[11pt,a4paper]{amsart}
\usepackage[T1]{fontenc}
\usepackage[utf8]{inputenc}
\usepackage{lmodern}
\usepackage{amsmath,amssymb,amsthm,mathtools}
\usepackage[margin=28mm]{geometry}
\usepackage{microtype}
\usepackage[colorlinks=true,linkcolor=blue,citecolor=blue,urlcolor=blue]{hyperref}
\newtheorem{theorem}{Theorem}
\newtheorem{lemma}[theorem]{Lemma}
\newtheorem{proposition}[theorem]{Proposition}
\theoremstyle{remark}
\newtheorem{remark}[theorem]{Remark}
\newcommand{\C}{\mathbb C}
\newcommand{\N}{\mathbb N}
\newcommand{\B}{\mathbb B}
\newcommand{\A}{\mathcal A}
\newcommand{\R}{\mathcal R}
\newcommand{\dist}{\operatorname{dist}}

\title[An asymptotic ball algebra]{An asymptotic ball algebra with countably determined spectra and a nonseparable character space}
\author{Juliusz Pham}
\date{November 2, 2026}
\begin{document}
\begin{abstract}
Let $X=\partial\B_2\subset\C^2$ and let $\R$ be the restrictions to $X$ of functions in the ball algebra $A(\overline\B_2)$. The unital uniform Banach algebra
\[
\A=\{(f_n)\in\ell^\infty(\N,C(X)):\dist_{C(X)}(f_n,\R)\longrightarrow0\}
\]
is inverse-closed in $\ell^\infty(\N,C(X))$. Consequently one countable family of evaluation characters generates the spectrum of every element of $\A$. Nevertheless $\Delta(\A)$ contains continuum many pairwise disjoint nonempty open subsets and is therefore nonseparable. This gives a negative answer to a problem posed by Ohrysko, Wojciechowski and Graham.
\end{abstract}
\maketitle

\section{The statement and the base algebra}

For a commutative unital complex Banach algebra $A$, let $\Delta(A)$ denote its character space with the Gelfand topology. A subset $S\subset\Delta(A)$ generates all spectra if
\[
\sigma_A(a)=\overline{\{\varphi(a):\varphi\in S\}}
\qquad(a\in A).
\]
Problem 18 of \cite{OWG} asks whether the existence of such a countable $S$ forces $\Delta(A)$ to be separable.

In Nikolski's terminology \cite[Definition~0.2.1]{Nikolski}, the spectrum of $A$ is $1$-visible from $S\subset\Delta(A)$ if $\sigma_A(a)=\overline{\{\varphi(a):\varphi\in S\}}$ for every $a\in A$. The classical Wiener--Pitt phenomenon, observed by Wiener and Pitt \cite{WP} and proved rigorously by \v{S}re\u{\i}der \cite{Sreider}, illustrates the failure of this property for Fourier evaluations on the convolution algebra $M(\mathbb T)$ of complex Borel measures on the circle. There exist measures $\mu$ such that $\overline{\{\widehat\mu(n):n\in\mathbb Z\}}\subsetneq\sigma_{M(\mathbb T)}(\mu)$. Our construction shows that $1$-visibility from a countable family is nevertheless compatible with a nonseparable character space.

\begin{theorem}\label{thm:main}
There is a commutative unital semisimple uniform Banach algebra $\A$ such that all spectra are generated by one countable family of characters, but $\Delta(\A)$ has continuum many pairwise disjoint nonempty open subsets.
\end{theorem}

Write
\[
\B_2=\{z=(z_1,z_2)\in\C^2:|z_1|^2+|z_2|^2<1\},\qquad X=\partial\B_2.
\]
Let $A(\overline\B_2)$ be the algebra of continuous functions on $\overline\B_2$ that are holomorphic in $\B_2$, and put
\[
\R=\{F|_X:F\in A(\overline\B_2)\}\subset C(X).
\]
The maximum-modulus principle gives $\|F\|_{\overline\B_2}=\|F|_X\|_X$. Thus restriction is an isometric algebra isomorphism onto the closed unital subalgebra $\R$. For $h\in\R$ denote its unique ball-algebra extension by $\widetilde h$, and define
\[
\chi_0(h)=\widetilde h(0).
\]
This is a character of norm one. If $z_i\in\R$ denotes the $i$th coordinate function on $X$, then
\[
\chi_0(z_1)=\chi_0(z_2)=0,\qquad |z_1(x)|^2+|z_2(x)|^2=1\quad(x\in X).
\tag{1}\label{eq:sphere}
\]

\begin{lemma}\label{lem:fullR}
The subalgebra $\R$ is inverse-closed in $C(X)$.
\end{lemma}
\begin{proof}
Let $h=F|_X\in\R$ and assume $h$ has no zeros on $X$. By continuity and compactness there is $0<r<1$ such that $F$ has no zeros on the shell $r<\|z\|\le1$. The function $1/F$ is holomorphic on $\B_2\setminus r\overline\B_2$. The classical Hartogs extension theorem extends it to a holomorphic function $G$ on $\B_2$; see, for example, \cite[Theorem~2.3.2]{Hormander}. Since $FG=1$ on the shell, the identity theorem gives $FG=1$ throughout $\B_2$. The function $G$ agrees with $1/F$ near the boundary and hence extends continuously to $\overline\B_2$. Therefore $1/h=G|_X\in\R$.
\end{proof}

Notice that inverse-closedness here is a statement about single functions. The pair $z_1,z_2$ has no common zero on $X$, although $\chi_0$ annihilates both.

\section{The asymptotic algebra and its spectra}

Give $E=\ell^\infty(\N,C(X))$ its usual pointwise product and norm
\[
\|f\|_E=\sup_n\|f_n\|_X.
\]
Define
\[
\A=\{f=(f_n)\in E:\dist(f_n,\R)\longrightarrow0\}.
\tag{2}\label{eq:A}
\]
Equivalently,
\[
\A=\ell^\infty(\N,\R)+c_0(\N,C(X)).
\tag{3}\label{eq:sum}
\]
Indeed, approximants $h_n\in\R$ with $\|f_n-h_n\|_X\to0$ can be chosen bounded whenever $f$ is bounded.

\begin{lemma}\label{lem:Aclosed}
The set $\A$ is a norm-closed unital subalgebra of $E$ and is inverse-closed in $E$.
\end{lemma}
\begin{proof}
The set $c_0(\N,C(X))$ is an ideal of $E$, and $\ell^\infty(\N,\R)$ is a unital subalgebra. Hence \eqref{eq:sum} is a unital subalgebra. To check closedness, suppose $f^{(k)}\to f$ uniformly in $E$. The inequality
\[
\dist(f_n,\R)\le\|f_n-f_n^{(k)}\|_X+\dist(f_n^{(k)},\R)
\]
gives $\limsup_n\dist(f_n,\R)\le\|f-f^{(k)}\|_E$. Letting $k\to\infty$ proves $f\in\A$.

For inverse-closedness, let $f\in\A$ be invertible in $E$, so that $|f_n(x)|\ge\delta>0$ for all $n$ and $x$. Choose bounded $h_n\in\R$ such that $\epsilon_n=\|f_n-h_n\|_X\to0$. For all sufficiently large $n$, $|h_n|\ge\delta/2$ on $X$. Lemma~\ref{lem:fullR} gives $h_n^{-1}\in\R$, and
\[
\|f_n^{-1}-h_n^{-1}\|_X\le\frac{2\epsilon_n}{\delta^2}\longrightarrow0.
\]
The finitely many initial coordinates do not affect condition \eqref{eq:A}. Thus $f^{-1}\in\A$.
\end{proof}

For $n\in\N$ and $x\in X$ let $\varepsilon_{n,x}(f)=f_n(x)$. These are characters and separate the elements of $\A$, so $\A$ is semisimple. Its norm is uniform, since it is a closed subalgebra of the commutative $C^*$-algebra $E$.

\begin{proposition}\label{prop:spectra}
Let $D\subset X$ be any countable dense subset. Then
\[
S=\{\varepsilon_{n,d}:n\in\N,\ d\in D\}
\]
is countable and generates all spectra in $\A$. More explicitly,
\[
\sigma_{\A}(f)=\overline{\bigcup_{n\in\N}f_n(X)}
=\overline{\{f_n(d):n\in\N,\ d\in D\}}.
\tag{4}\label{eq:spectrum}
\]
\end{proposition}
\begin{proof}
Each value $f_n(x)$ belongs to the spectrum by evaluation at $\varepsilon_{n,x}$, so the closure of their union is contained in $\sigma_{\A}(f)$. If $\lambda$ is outside that closure, then $|f_n(x)-\lambda|$ has a uniform positive lower bound. Lemma~\ref{lem:Aclosed} shows that $(f-\lambda 1)^{-1}\in\A$. This proves the first equality. The second follows from continuity of every $f_n$ and density of $D$.
\end{proof}

Restriction to the countable set $\N\times D$ embeds $E$, and hence $\A$, isometrically into $\ell^\infty(\N\times D,\C)$. The image of $\A$ is closed and inverse-closed there, by the same uniform lower-bound argument. Since $\N\times D$ is countable, this proves the last realization assertion of Theorem~\ref{thm:main}.

\section{Finite-coordinate and ultrafilter characters}

For $L\subset\N$ put
\[
(e_L)_n=1_L(n)\,1_X.
\]
Every $e_L$ belongs to $\ell^\infty(\N,\R)\subset\A$ and is an idempotent. Write $e_n=e_{\{n\}}$, and let $Z_i\in\A$ be the constant sequence whose every term is $z_i$, for $i=1,2$.

\begin{lemma}\label{lem:finite}
If $\varphi\in\Delta(\A)$ and $\varphi(e_n)=1$, then $\varphi=\varepsilon_{n,x}$ for some $x\in X$. In particular,
\[
|\varphi(Z_1)|^2+|\varphi(Z_2)|^2=1.
\]
\end{lemma}
\begin{proof}
The algebra $e_n\A$ is exactly $C(X)$ supported in coordinate $n$, since every such sequence belongs to $c_0(\N,C(X))$. Restricting $\varphi$ gives a unital character of $C(X)$, hence evaluation at some $x$. For every $f\in\A$, multiplicativity gives $\varphi(f)=\varphi(e_nf)=f_n(x)$. Now apply \eqref{eq:sphere}.
\end{proof}

\begin{lemma}\label{lem:ultra}
For every free ultrafilter $p$ on $\N$ there is a character $\Phi_p\in\Delta(\A)$ satisfying
\[
\Phi_p(Z_1)=\Phi_p(Z_2)=0,\qquad
\Phi_p(e_L)=\begin{cases}1,&L\in p,\\0,&L\notin p.\end{cases}
\tag{5}\label{eq:ultra}
\]
\end{lemma}
\begin{proof}
For $f\in\A$, choose bounded $h_n\in\R$ with $\|f_n-h_n\|_X\to0$ and set
\[
\Phi_p(f)=\lim_{n\to p}\chi_0(h_n).
\tag{6}\label{eq:defultra}
\]
If $k_n$ is a second choice, then
\[
|\chi_0(h_n)-\chi_0(k_n)|\le\|h_n-k_n\|_X\longrightarrow0.
\]
A free ultrafilter contains every cofinite set, so \eqref{eq:defultra} is independent of the choice. Linearity follows by using $h_n+k_n$ as approximants to a sum. For products, if $h_n$ approximates $f_n$ and $k_n$ approximates $g_n$, then $h_nk_n$ approximates $f_ng_n$, because the sequences are bounded. The multiplicativity of $\chi_0$ and of scalar ultrafilter limits gives multiplicativity of $\Phi_p$. Also $\Phi_p(1)=1$ and $|\Phi_p(f)|\le\|f\|_E$. Finally, choose the exact approximants $z_i$ and $1_L(n)1_X$ to obtain \eqref{eq:ultra}.
\end{proof}

\section{Continuum many disjoint open subsets}

Set
\[
V=\{\varphi\in\Delta(\A):|\varphi(Z_1)|^2+|\varphi(Z_2)|^2<1/4\}.
\]
This is open and nonempty by Lemma~\ref{lem:ultra}. Lemma~\ref{lem:finite} shows that every $\varphi\in V$ annihilates every $e_n$, and therefore every finite sum of the $e_n$.

Choose a family $\mathcal F$ of continuum many infinite subsets of $\N$ with pairwise finite intersections. Such a family is obtained by enumerating the finite binary words and taking the sets of prefixes of distinct infinite binary sequences. For $L\in\mathcal F$ define
\[
U_L=V\cap\{\varphi:|\varphi(e_L)-1|<1/2\}.
\tag{7}\label{eq:UL}
\]
These sets are open. Since $e_L$ is idempotent, $\varphi(e_L)\in\{0,1\}$; hence the second condition in \eqref{eq:UL} simply means $\varphi(e_L)=1$.

Each $U_L$ is nonempty: choose a free ultrafilter $p$ containing $L$ and use $\Phi_p$. If $L\ne M$ and $\varphi\in U_L\cap U_M$, multiplicativity gives
\[
1=\varphi(e_L)\varphi(e_M)=\varphi(e_{L\cap M}).
\]
But $L\cap M$ is finite, and every character in $V$ annihilates finite-coordinate idempotents. This is a contradiction. The sets $U_L$ are therefore pairwise disjoint.

A countable dense subset would have to meet each of these continuum many pairwise disjoint nonempty open sets. Consequently $\Delta(\A)$ is nonseparable, completing the proof of Theorem~\ref{thm:main}.

\begin{remark}[Two-function form of the open sets]
For $L\subset\N$ put
\[
A_L=1-e_L+e_LZ_1,\qquad B_L=e_LZ_2.
\]
At every coordinate evaluation one has
\[
|\varepsilon_{n,x}(A_L)|^2+|\varepsilon_{n,x}(B_L)|^2=1.
\]
Nevertheless, if $p$ is free and $L\in p$, then $\Phi_p(A_L)=\Phi_p(B_L)=0$. Moreover,
\[
U_L=\{\varphi:|\varphi(A_L)|^2+|\varphi(B_L)|^2<1/4\}.
\]
Thus the disjoint open sets require two simultaneous conditions. There is no conflict with \eqref{eq:spectrum}, which concerns one element at a time.
\end{remark}

\begin{remark}[The Shilov boundary]
The Shilov boundary of $\A$ is exactly the closure of $\{\varepsilon_{n,x}:n\in\N,x\in X\}$, and is separable. Indeed, this closure is a closed norming set. Conversely, finite-coordinate functions from $C(X)$ provide functions supported in arbitrarily small neighborhoods of each $\varepsilon_{n,x}$, showing that all these evaluations are strong boundary points. The open set $V$ lies outside this closed boundary. In particular, nonseparability is not being inferred merely from a nonseparable closed subset of a possibly separable space: it follows from the explicit disjoint open sets above.
\end{remark}

\begin{remark}[A general construction principle]
The same proof works for any compact Hausdorff separable $X$ and any closed unital inverse-closed subalgebra $R\subset C(X)$ admitting a character $\chi$ outside the evaluation image of $X$. By compactness, there are finitely many $r_i\in R$ with $\chi(r_i)=0$ and $\inf_{x\in X}\sum_i|r_i(x)|^2>0$. Replace $z_1,z_2$ by these $r_i$ throughout, and use a countable dense subset of $X$ in Proposition~\ref{prop:spectra}. No metrizability is used. The ball algebra on $\partial\B_2$ supplies an explicit instance of all the hypotheses.
\end{remark}

\begin{remark}[The global implication for corona theorems]
The preceding amplification principle supplies the additional argument needed to connect Problem 18 with corona theorems. Let $\Omega$ be a separable complex domain and let $H=H^\infty(\Omega)$. Put $F=\overline{\{\delta_z:z\in\Omega\}}\subset\Delta(H)$. The space $F$ is compact and separable: evaluations at a countable dense subset of $\Omega$ are dense in it. Restriction of the Gelfand transform identifies $H$ isometrically with a closed unital subalgebra $R\subset C(F)$. This subalgebra is inverse-closed: a uniform lower bound on $F$ is a uniform lower bound on $\Omega$, so the reciprocal is again bounded and holomorphic. If the corona is nonempty, then $\Delta(R)$ has a character outside $F$. The amplification principle therefore produces an algebra with countably generated spectra and nonseparable character space. Consequently, a positive answer to Problem 18 for all commutative unital Banach algebras would imply the corona theorem for $H^\infty(\Omega)$, including balls and polydisks. This implication is stronger than merely applying a putative positive answer directly to $H^\infty(\Omega)$, which would give only separability of its character space.
\end{remark}

\section*{Use of artificial intelligence}
The proof was obtained in collaboration with OpenAI’s GPT 6.1 Sol model. The author has read the entire manuscript and takes full responsibility for its contents.

\end{document}